\documentclass[11pt]{article}

\usepackage{amssymb}
\usepackage[english]{babel}
\usepackage{graphicx}

\newtheorem{theo}{Theorem}[section]
\newtheorem{propo}[theo]{Proposition}
\newtheorem{lema}[theo]{Lemma}

\newtheorem{coro}[theo]{Corollary}

\def\proof{{\boldmath  $Proof.$}\hskip 0.3truecm}
\def\endproof{\quad $\Box$}

\input amssym.def
\newsymbol\rtimes 226F
\newfont{\nset}{msbm10}

\def\A{{\mbox {\boldmath $A$}}}
\def\B{{\mbox {\boldmath $B$}}}
\def\C{{\mbox {\boldmath $C$}}}

\def\E{{\mbox {\boldmath $E$}}}

\def\I{{\mbox {\boldmath $I$}}}
\def\J{{\mbox {\boldmath $J$}}}

\def\A{{\mbox {\boldmath $A$}}}

\def\matrix0{{\mbox {\boldmath $O$}}}
\def\M{\mbox{\boldmath $M$}}

\def\e{{\mbox{\boldmath $e$}}}

\def\ubar{\bar{u}}

\def\vbar{\bar{v}}
\def\wbar{\bar{w}}
\def\zbar{\bar{z}}
\def\x{{\mbox{\boldmath $x$}}}
\def\y{{\mbox{\boldmath $y$}}}

\def\z{{\mbox{\boldmath $z$}}}
\def\vec0{\mbox{\bf 0}}

\def\dist{\mathop{\partial }\nolimits}

\def\Ker{\mathop{\rm Ker }\nolimits}

\def\sp{\mathop{\rm sp }\nolimits}

\begin{document}

\title{On Perturbations of Almost Distance-Regular Graphs\footnote{This version is published in Linear Algebra and its Applications 435 (2011), 2626-2638.}}
\author{C. Dalf\'{o}$^\dag$, E.R. van Dam$^\ddag$ and M.A. Fiol$^\dag$
\\ \\
{\small $^\dag$Universitat Polit\`ecnica de Catalunya, Dept. de Matem\`atica Aplicada IV} \\
{\small Barcelona, Catalonia} {\small (e-mails: {\tt
\{cdalfo,fiol\}@ma4.upc.edu})} \\
{\small $^\ddag$Tilburg University, Dept. Econometrics and O.R.} \\
{\small Tilburg, The Netherlands} {\small (e-mail: {\tt
edwin.vandam@uvt.nl})} \\
}
\date{In honour of Drago\v{s} Cvetkovi\'c, on his 70th birthday}
\maketitle

\noindent {\footnotesize Keywords: Distance-regular graph, Walk-regular graph,
Eigenvalues, Perturbation, Cospectral graphs

}

\begin{abstract}
In this paper we show that certain almost distance-regular
graphs, the so-called $h$-punctually walk-regular graphs, can
be characterized through the cospectrality of their perturbed
graphs. A graph $G$ with diameter $D$ is called $h$-punctually
walk-regular, for a given $h\le D$, if the number of paths of
length $\ell$ between a pair of vertices $u,v$ at distance $h$
depends only on $\ell$. The graph perturbations considered here
are deleting a vertex, adding a loop, adding a pendant edge,
adding/removing an edge, amalgamating vertices, and adding a
bridging vertex. We show that for walk-regular graphs some of
these operations are equivalent, in the sense that one
perturbation produces cospectral graphs if and only if the
others do. Our study is based on the theory of graph
perturbations developed by Cvetkovi\'c, Godsil, McKay,
Rowlinson, Schwenk, and others. As a consequence, some new
characterizations of distance-regular graphs are obtained.
\end{abstract}

\section{Introduction}

Both the theory of distance-regular graphs and that of graph perturbations have
been widely developed in the last decades. The importance of the former can be
grasped from the comment in the preface of the comprehensive monograph of
Brouwer, Cohen, and Neumaier \cite{bcn}: ``Most finite objects bearing `enough
regularity' are closely related to certain distance-regular graphs.'' Thus,
many characterizations of a combinatorial and algebraic nature of
distance-regular graphs are known (see \cite{f02}), and they have given rise to
several generalizations, such as association schemes (see Brouwer and Haemers
\cite{bh95}) and almost distance-regular graphs \cite{ddfgg10}. With respect to
the latter, the spectral properties of modified (or `perturbed')  graphs have
relevance in Chemistry, in the construction of isospectral molecules, as well
as in other areas of  graph  theory (as in the reconstruction conjecture); see
Cvetkovi\'c, Doob, and Sachs \cite{cds80}, Rowlinson \cite{r91per,r96}, and
Schwenk \cite{sch79}. The aim of this paper is to put together different ideas
and results from both theories to show that certain almost distance-regular
graphs, the so-called $h$-punctually walk-regular (or $h$-punctually
spectrum-regular) graphs, can be characterized through the cospectrality of
their perturbed graphs. We consider three one-vertex perturbations, namely,
vertex deletion, adding a loop at a vertex, and adding a pendant edge at a
vertex. These three perturbations are extended to pairs of vertices to obtain
two-vertex `separate' perturbations. We also consider three two-vertex `joint'
perturbations, namely adding/removing an edge, amalgamating two vertices, and
adding a bridging vertex. We show that for walk-regular graphs all these
two-vertex operations are equivalent, in the sense that one perturbation
produces cospectral graphs if and only if the others do. We also consider
perturbations on a set of vertices, and their impact on almost distance-regular
graphs. As a consequence, we obtain some new characterizations of
distance-regular graphs, in terms of the cospectrality of their perturbed
graphs.

\section{Preliminaries}

In this section we give the basic definitions, notation and
results on which our study is based. For completeness, we prove
again some known results. Accordingly, we also recall some
basic results on the computation of determinants which are used
in our study.

\subsection{Graphs and their spectra} \label{secgraphspectra}

Let $G=(V,E)$ be a (connected) graph with vertex set $V$ and edge set $E$. The
adjacency between vertices $u,v\in V$, that is $uv\in E$, is denoted by $u\sim
v$, and their distance is $\dist(u,v)$. Let $\A=(a_{uv})$ be the adjacency
matrix of $G$, with characteristic polynomial $\phi_G(x)$, and spectrum $\sp G
= \{\lambda_0^{m_0},\lambda_1^{m_1},\dots, \lambda_d^{m_d}\}$, where the
different eigenvalues of $G$ are in decreasing order,
$\lambda_0>\lambda_1>\cdots >\lambda_d$, and the superscripts stand for their
multiplicities $m_i=m(\lambda_i)$. For $i=0,1,\ldots, d$, let $\E_i$ be the
principal idempotent of $\A$, which corresponds to the orthogonal projection
onto the eigenspace  ${\cal E}_i=\Ker (\lambda_i\I-\A)$. In particular, if $G$
is regular, $\E_0=\frac{1}{n}\J$, where $\J$ stands for the all-$1$ matrix. As
is well known, the idempotents satisfy the following properties:
$\E_i\E_j=\delta_{ij}\E_i$ (with $\delta_{ij}$ being the Kronecker delta),
$\A\E_i=\lambda_i\E_i$, and $q(\A)=\sum_{i=0}^d q(\lambda_i)\E_i$ for every
rational function $q$ that is well-defined at each eigenvalue of $\A$; see, for
instance, Godsil \cite{g93}. The $uv$-entry $m_{uv}(\lambda_i)=(\E_i)_{uv}$ of
the idempotent $\E_i$ is called the {\em crossed $(uv$-$)$local multiplicity}
of $\lambda_i$. As some direct consequences of the above properties, the
following lemma gives some properties of these parameters (see, for example,
\cite{f01}).
\begin{lema}
\label{crossed-mul} For $u,v\in V$, the crossed local
multiplicities of each eigenvalue $\lambda_i$,
$i=0,1,\ldots,d$, satisfy the following properties:
\begin{itemize}
\item[$(a)$] $\sum_{i=0}^d m_{uv}(\lambda_i)=\delta_{uv}$.
\item[$(b)$] $\sum_{w\sim v} m_{uw}(\lambda_i)=\lambda_i m_{uv}(\lambda_i)$.
\item[$(c)$] $a_{uv}^{(\ell)}=(\A^{\ell})_{uv}=\sum_{i=0}^d m_{uv}(\lambda_i)\lambda_i^{\ell}$.
\end{itemize}
\end{lema}
Note that the $uv$-entry $a_{uv}^{(\ell)}$ of the power matrix $\A^{\ell}$ is
equal to the number of walks of length $\ell$ between vertices $u,v$. Rowlinson
\cite{r97} showed that a graph $G$ is distance-regular if and only if this
number of walks only depends on $\ell=0,1,\ldots,d$ and the distance
$\dist(u,v)$ between $u$ and $v$. Similarly, $G$ is distance-regular if and
only if its local crossed multiplicities $m_{uv}(\lambda_i)$ only depend on
$\lambda_i$ and $\dist(u,v)$; see \cite{f02}. Inspired by these
characterizations, the authors \cite{ddfgg10} introduced the following concepts
as different approaches to `almost distance-regularity'. We say that a graph
$G$ with diameter $D$ and $d+1$ distinct eigenvalues is $h$-\emph{punctually
walk-regular}, for a given $h\le D$, if for every $\ell \ge 0$ the number of
walks of length $\ell$ between a pair of vertices $u,v$ at distance
$\dist(u,v)=h$ does not depend on $u,v$. Similarly, we say  that $G$ is
\emph{$h$-punctually spectrum-regular}, for a given $h\leq D$ if for all $i \le
d$, the crossed $uv$-local multiplicities of $\lambda_i$ are the same for all
pairs of vertices $u,v$ at distance $\dist(u,v)=h$. In this case, we write
$m_{uv}(\lambda_i)=m_{hi}$. The concepts of $h$-punctual walk-regularity and
$h$-punctual spectrum-regularity are equivalent. For $h=0$, the concepts are
equivalent to walk-regularity (a concept introduced by Godsil and McKay in
\cite{gmk80}) and spectrum-regularity (see Fiol and Garriga \cite{fg98}),
respectively.

\subsection{Graph perturbations}

As mentioned above, we consider three basic graph
perturbations which involve a given vertex $u\in V$:
\begin{itemize}
\item[{\bf P1.}] $G-u$ is the graph obtained from $G$ by
    removing $u$ and all the edges incident to it.
\item[{\bf P2.}] $G+uu$ is the (pseudo)graph obtained from
    $G$ by adding a loop at $u$. (In this case the graph
    obtained has adjacency matrix as expected, with its
    $uu$-entry equal to $1$.)
\item[{\bf P3.}] $G+u\ubar$ is the graph obtained from $G$
    by adding a pendant edge at $u$ (thus creating a new
    vertex $\ubar$).
\end{itemize}
Two vertices $u,v$ satisfying $\sp (G-u)=\sp(G-v)$ were called {\em cospectral}
by Herndon and Ellzey \cite{he75}. We say that a graph is {\em $0$-punctually
cospectral} when all its vertices are cospectral; a concept that we will
generalize below. It is well-known that a graph is $0$-punctually cospectral if
and only if it is walk-regular; see Proposition \ref{propo0}, where we also
relate this to the perturbations {\bf P2} and {\bf P3}. In fact, the proof of
Proposition \ref{propo0} implies that cospectral vertices $u,v$ can be
equivalently defined by requiring that $\sp (G+uu)=\sp(G+vv)$ or $\sp
(G+u\ubar)=\sp(G+v\vbar)$.

Given a vertex subset $U\subset V$, we can also consider the graphs obtained by
applying any of the above perturbations to every vertex of $U$, with natural
notation $G-U$, $G+UU$ and $G+U\bar{U}$. In particular, when $U=\{u,v\}$, we
also write $G-u-v$, $G+uu+vv$ and $G+u\ubar+v\vbar$.

Building on the concept of cospectral vertices, Schwenk
\cite{sch79} considered the analogue for sets: Two vertex
subsets $U,U'\subset V$ are  {\em removal-cospectral} if there
exists a one-to-one mapping $U \rightarrow U'$ such that, for
every $W\subset U$, the graphs $G-W$ and $G-W'$ are cospectral. A
main result of his paper was the following necessary condition
for two sets being removal-cospectral:

\begin{theo} \label{thmschwenk}
\cite{sch79} If $U,U'$ are removal-cospectral sets, then
$a_{uv}^{(\ell)}=a_{u'v'}^{(\ell)}$  for all pairs of vertices
$u,v\in U$ and all $\ell\ge 0$.
\end{theo}
Godsil \cite{g92} proved that two vertex subsets $U,U'$ are
removal-cospectral if and only if for every subset $W\subset U$
with at most two vertices, the subsets $W,W'$ are
removal-cospectral (for both an alternative proof and a
geometric interpretation of this result, see Rowlinson
\cite{r96}).

As a consequence of Theorem \ref{thmschwenk}, notice that for
$\{u,v\}$ and $\{u',v'\}$ to be removal-cospectral we need that
$\dist(u,v)=\dist(u',v')$. Otherwise, if
$r=\dist(u,v)<\dist(u',v')$, say, we would have
$a_{uv}^{(r)}>0$ whereas $a_{u'v'}^{(r)}=0$. Inspired by this
property, we say that two vertex subsets are {\em isometric}
when there exists a one-to-one mapping $U \rightarrow U'$ such that, for every pair $u,v\in U$, we have $\dist(u,v)=\dist(u',v')$.
So, if two sets are removal-cospectral then they are also
isometric. In the last section, we will show that the converse
is also true for distance-regular graphs.

\begin{figure}[h!]
\begin{center}
\resizebox{30mm}{!}{\includegraphics{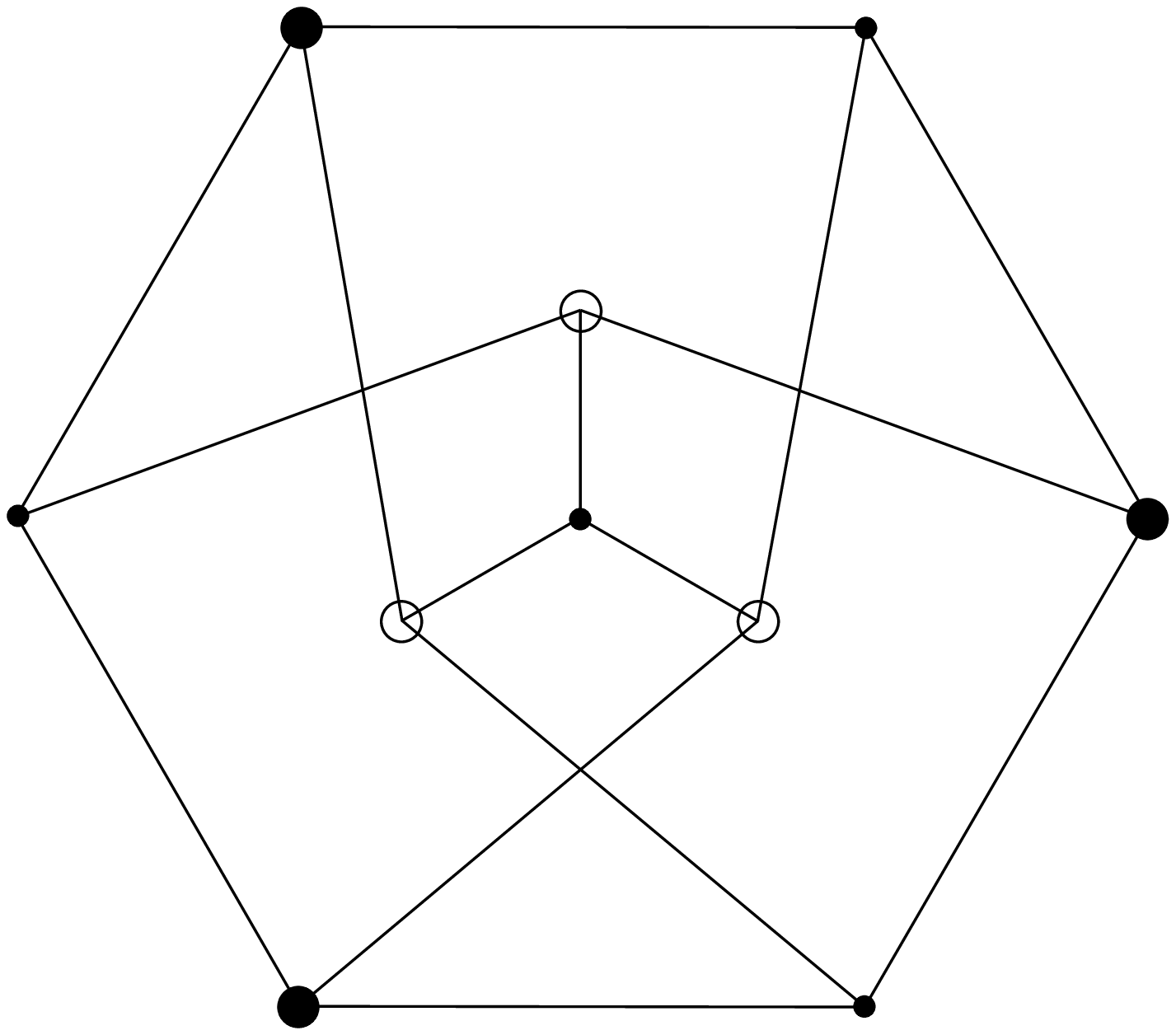}} \caption{Petersen graph with 3-cocliques} \label{Petersen}
\end{center}
\end{figure}

\begin{figure}[h!]
\begin{center}
\resizebox{30mm}{!}{\includegraphics{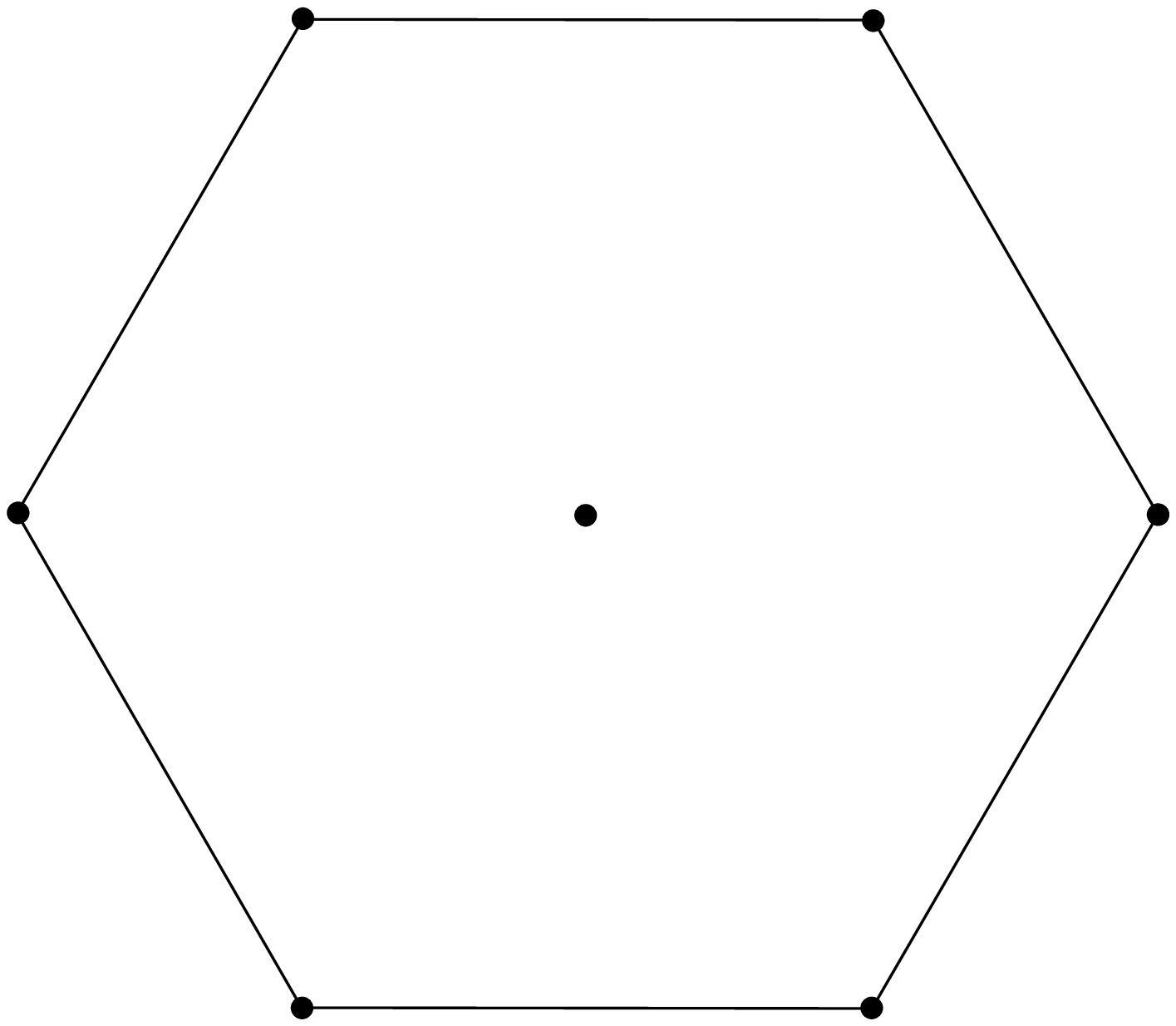}} \nobreak \hspace{0.5 cm}
\resizebox{22mm}{!}{\includegraphics{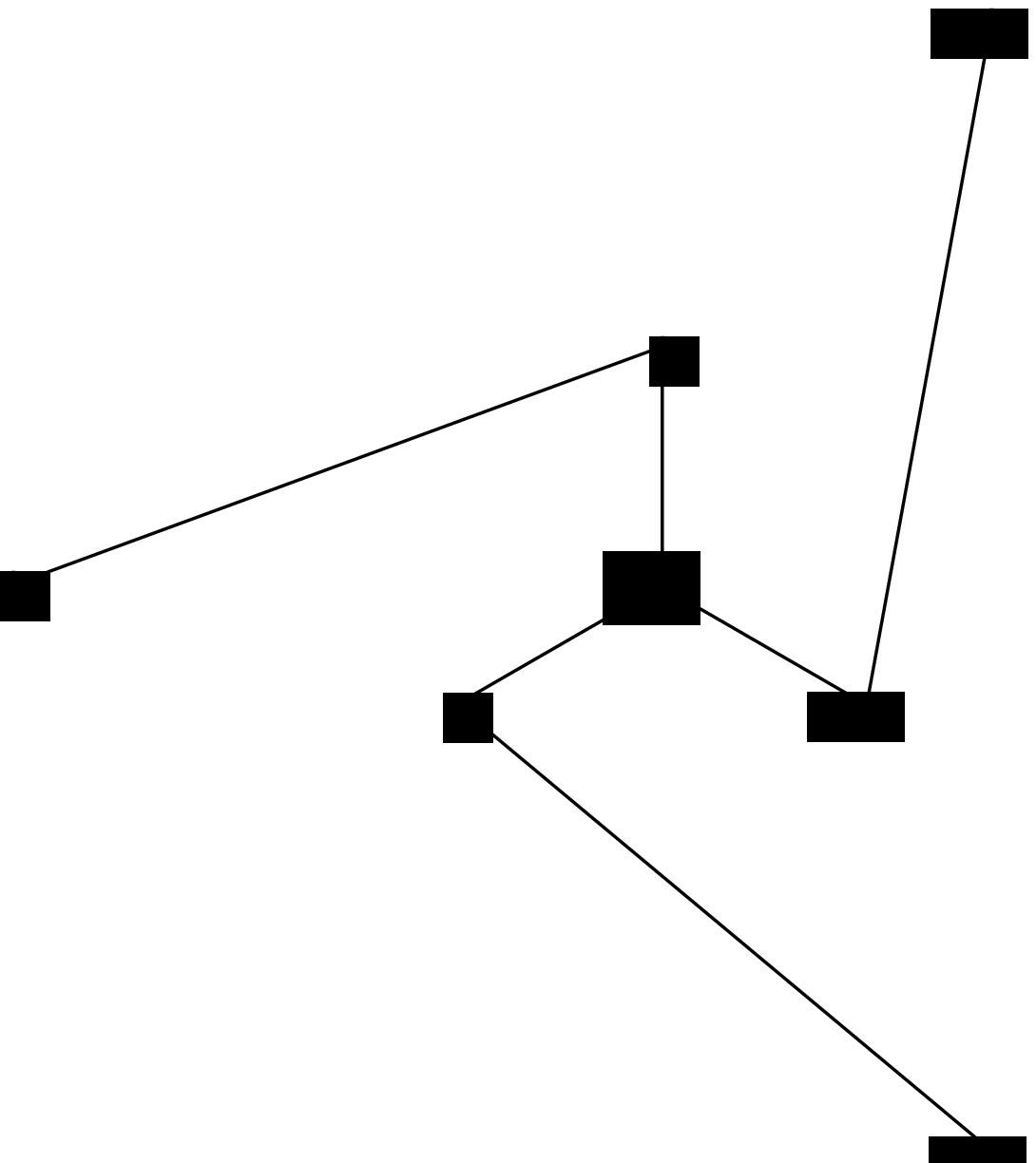}} \nobreak \hspace{2 cm}
\resizebox{30mm}{!}{\includegraphics{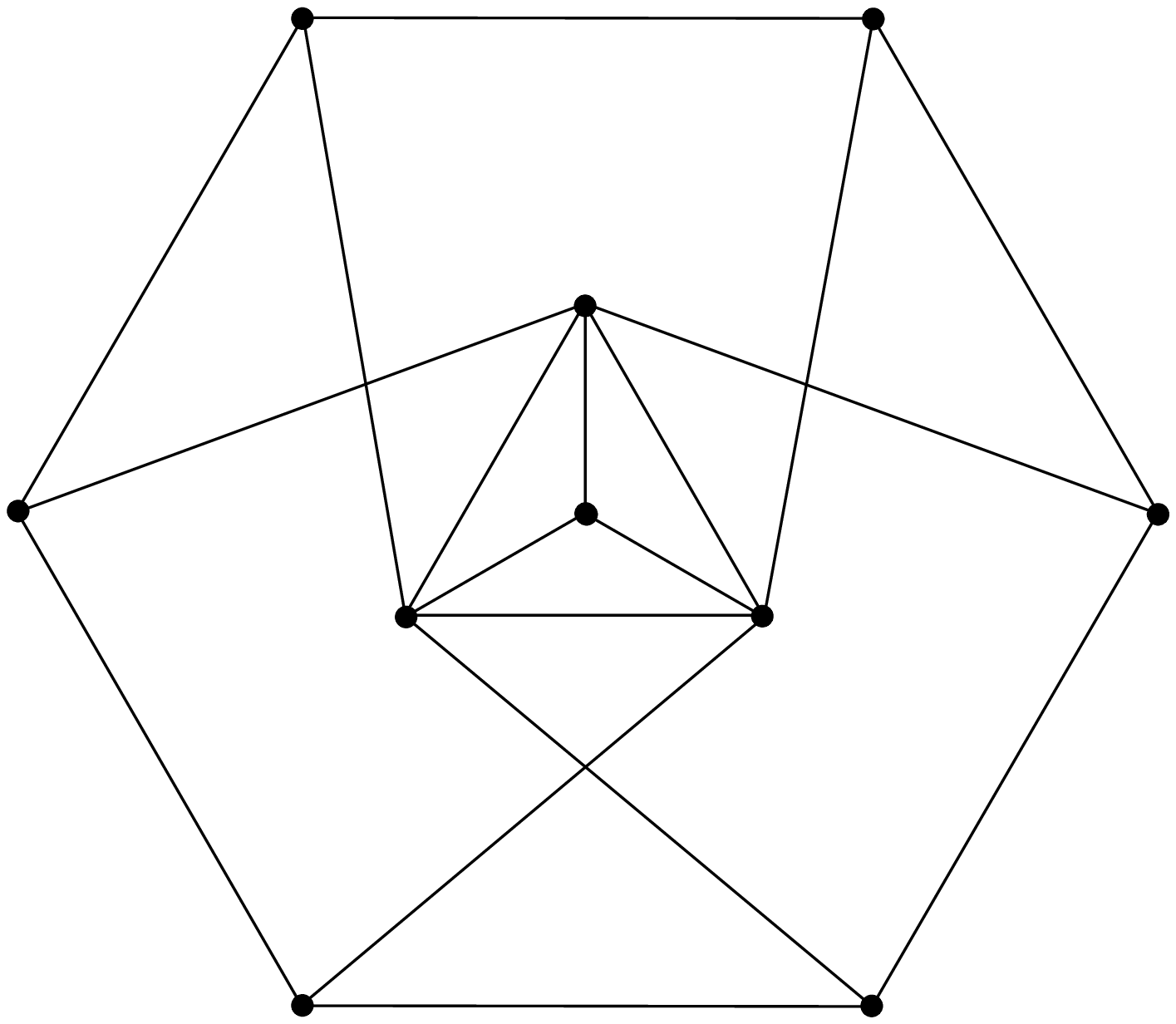}} \nobreak \hspace{0.5 cm}
\resizebox{30mm}{!}{\includegraphics{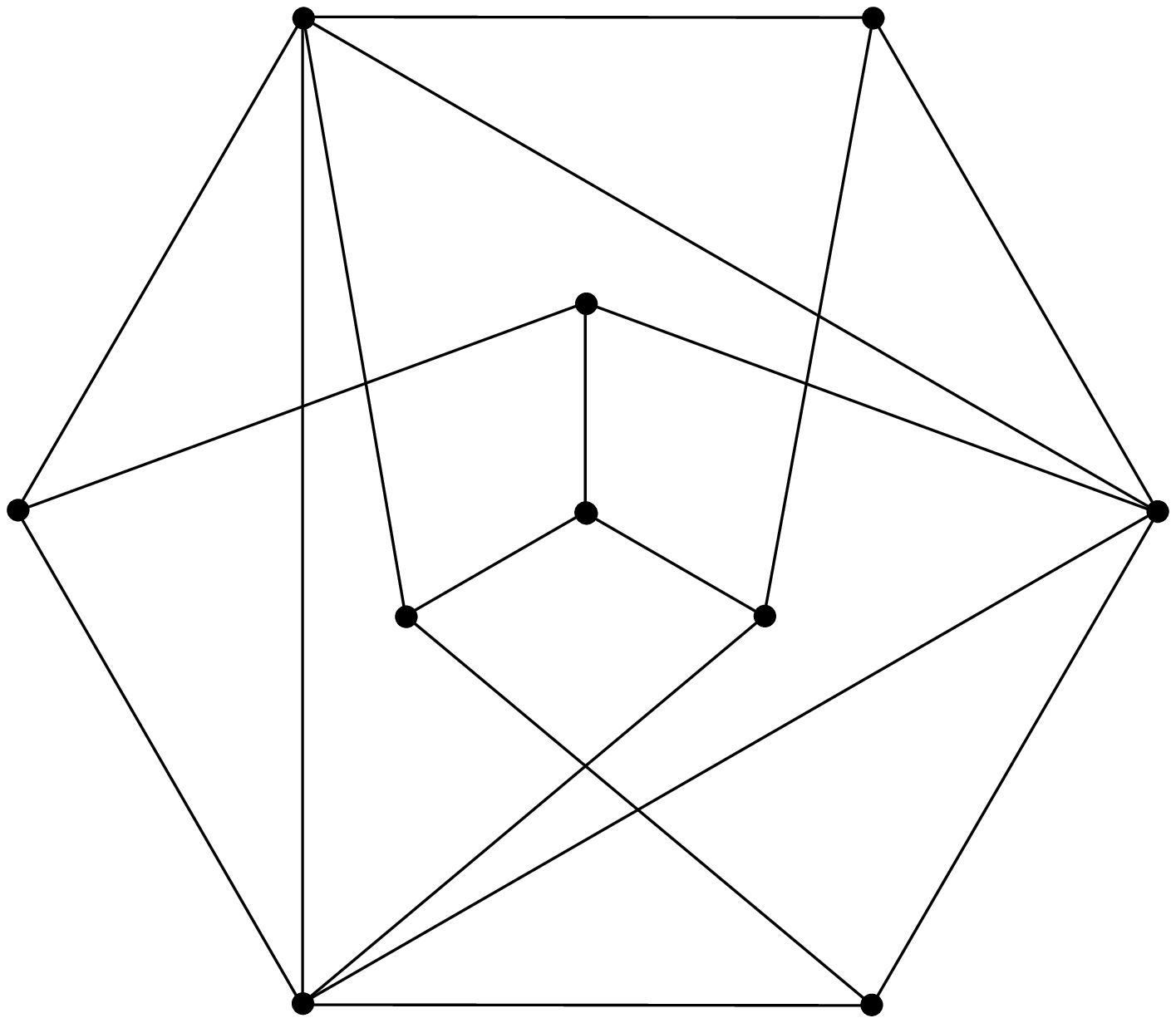}}
\end{center}
\caption{Two pairs of cospectral graphs: removing vertices and adding edges}
\label{cocliquesremoved}
\end{figure}

For example, in the Petersen graph all cocliques (that is, independent sets) of
size $3$ are removal-cospectral. Since there are two different kinds of such
cocliques (one of these is indicated in Figure \ref{Petersen} by the empty
dots, and the other by the thick dots), removing them gives a pair of
cospectral but non-isomorphic graphs. This is the left pair in Figure
\ref{cocliquesremoved}. The right pair is obtained by adding edges to the
cocliques. This also gives cospectral but non-isomorphic graphs since, as was
proved by Schwenk \cite{sch79}, if $U$ and $U'$ are removal-cospectral sets,
then any graph may be attached to all the points of $U$ and to the points of
$U'$ with the two graphs so formed being cospectral.

In our framework of almost distance-regular graphs, the case when the two
vertices of $W$ are at a given distance proves to be specially relevant, and
leads us to the following definition: A graph $G$ with diameter $D$ is {\em
$h$-punctually cospectral}, for a given $h\le D$, when, for {\em all} pairs of
vertices $u,v$ and $w,z$, both at distance $\dist(u,v)= \dist(w,z)=h$, we have
$\sp(G-u-v)=\sp(G-w-z)$. Again, we will show later (in Lemma \ref{co-sp}) that
this concept can also be defined by using the other graph perturbations
considered here. Notice that, since there are no restrictions on either pair of
vertices, except for their distance, this is equivalent to the sets $W=\{u,v\}$
and $W'=\{u',v'\}$, with both mappings $u'=w$, $v'=z$ and $u'=z$, $v'=w$, being
removal-cospectral.

Then, using our terminology, Schwenk's theorem implies the
following corollary:
\begin{coro}
\label{coroSchwenk} If a graph $G$ is $j$-punctually cospectral for $j=0,h$,
then it is $j$-punctually walk-regular for $j=0,h$.
\end{coro}
Answering a question of Schwenk \cite{sch79}, Rowlinson
\cite{r96} proved the following characterization of
removal-cospectral sets, which we give in terms of the local
crossed multiplicities:

\begin{theo}
\label{theo-row} \cite{r96} The vertex (non-empty) subsets
$U,U'$ are removal-cospectral if and only if
$m_{uv}(\lambda_i)=m_{u'v'}(\lambda_i)$  for all $u,v\in U$ and
$i=0,1,\ldots,d$.
\end{theo}
In fact, Rowlinson gave his result in terms of the so-called {\em star
sequences}  $\{\E_i\e_u: u\in U\}$, $i=0,1,\ldots, d$ and $\{\E_i\e_{u'}: u'\in
U'\}$, $i=0,1,\ldots, d$, where $\e_u$ stands for the $u$-th unit vector.

Again, in our context we have the following consequence:

\begin{coro} \label{coroequiv}
A graph $G$ is $j$-punctually cospectral for $j=0,h$ if and
only if it is $j$-punctually spectrum-regular for $j=0,h$.
\end{coro}
We remind the reader that the concepts of $h$-punctually
walk-regularity and $h$-punctually spectrum-regularity are
equivalent, so Corollary \ref{coroequiv} implies Corollary
\ref{coroSchwenk}. As this corollary is one of the crucial
characterizations for us, we will restate (and prove) it later on
as Theorem \ref{theo-punct}.

\subsection{Computing determinants}

Our first study will use the two following lemmas to compute
determinants. A proof of the first result can be found, for
instance, in Godsil \cite[p. 19]{g93}. For the argument for
Jacobi's determinant identity, see Rowlinson \cite[p. 212]{r96}, for example.

\begin{lema}
\label{det1} Let $\A$ and $\B$ be two $n\times n$ matrices.
Then, $\det(\A+\B)$ equals the sum of the determinants of the
$2^n$ matrices obtained by replacing every subset of the
columns of $\A$ by the corresponding subset of the columns of
$\B$.
\end{lema}
In particular, for all column vectors $\x$, $\y$ of size $n$
and $n\times (n-1)$ matrix $\M$, we have the well-known
linearity property
\begin{equation}
\label{column-linearity}
\det(\x+\y|\M)=\det(\x|\M)+\det(\y|\M).
\end{equation}

\begin{lema}
\label{det2} $($Jacobi's determinant identity$)$ Let $\A$ be an
invertible matrix with rows and columns indexed by the elements
of $V$. For a given nontrivial subset $U$ of $V$, let $\A[U]$ denote the
principal submatrix of $\A$ on $U$.  Let
$\overline{U}=V\setminus U$. Then,
$$
\det \A[U]=\det \A \det \A^{-1}[\overline{U}].
$$
\end{lema}

\section{Walk-regular graphs}

Our main results were inspired by the following
characterizations of walk-regular graphs:
\begin{propo}
\label{propo0}
The following statements are equivalent:
\begin{itemize}
\item[$(a)$]
$G$ is walk-regular
\item[$(b)$]
$G$ is spectrum-regular.
\item[$(c)$] $\sp (G-u) = \sp (G-v)$ for all vertices
    $u,v$.
\item[$(d)$] $\sp (G+uu) = \sp (G+vv)$ for all vertices
    $u,v$.
\item[$(e)$] $\sp (G+u\ubar) = \sp (G+v\vbar)$ for all
    vertices $u,v$.
\end{itemize}
\end{propo}
\proof Let $G$ have adjacency matrix $\A$. The equivalence
$(a)\iff (b)$ was proved by Delorme and Tillich \cite{dt97} and
Fiol and Garriga \cite{fg98}.

Godsil and McKay \cite{gmk81} obtained a relation between a walk-generating
function of $G$ and the characteristic polynomials of $G$ and $G-u$. This can
be formulated (see also \cite[p. 83]{CRS}) as
$$\phi_{G-u}(x)=\phi_G(x)\sum_{i=0}^d\frac{m_u(\lambda_i)}{x-\lambda_i}.$$
This can also be proved by using Lemma \ref{det2}. Indeed, let
$U=V\setminus\{u\}$ and $\C=x\I-\A$. Then, as
$\C^{-1}=\sum_{i=0}^d\frac{\E_i}{x-\lambda_i}$, we have
$$
\phi_{G-u}(x)  =  \det \C[U] =\det \C\det \C^{-1}[\{u\}]
   =  \phi_G(x) \sum_{i=0}^d \frac{\E_i[\{u\}]}{x-\lambda_i}
   =  \phi_G(x) \sum_{i=0}^d \frac{m_u(\lambda_i)}{x-\lambda_i}.
$$
Therefore $(b)\Rightarrow (c)$. Conversely, if
$\phi_{G-u}(x)=\phi_{G-v}(x)$, then the limit
$$
\lim_{x\rightarrow \lambda_i}\frac{\phi_{G-u}(x)}{\phi_{G-v}(x)}=\frac{m_u(\lambda_i)}{m_v(\lambda_i)}
$$
yields that $m_u(\lambda_i)=m_v(\lambda_i)$ for every
$i=0,1,\ldots,d$, hence $(c)\Rightarrow (b)$.

If we apply Lemma \ref{det1} to compute the determinant of
$(x\I-\A)+ (-\B)$, where $\B$ is the matrix with the only
non-zero entry $(\B)_{uu}=1$, we get
\begin{equation}
\label{G+uu vs G-u}
\phi_{G+uu}(x)=\phi_G(x)-\phi_{G-u}(x),
\end{equation}
thus proving the equivalence $(c)\iff (d)$.

Finally, the natural determinantal expansion of
$$
x\I-\A_{\ubar}=\left(
\begin{array}{cc}
x   &\mbox{$-1$\ \ \ \ } \vec0^{\top} \\
\begin{array}{c}
-1  \\
\vec0
\end{array}  & x\I-\A
\end{array}
\right),
$$
where $\A_{\ubar}$ is the adjacency matrix of $G+u\ubar$ (with
the first two rows and columns indexed by the vertices $\ubar$
and $u$), gives the well-known result
\begin{equation}
\label{G+uubar vs G-u}
\phi_{G+u\ubar}(x)=x\phi_G(x)-\phi_{G-u}(x)
\end{equation}
(see also, for instance, Rowlinson \cite{r91per}), thus proving
that $(c)\iff (e)$.
\endproof

Thus, we have just proved that a graph $G$ is ($0$-punctually) walk-regular or
($0$-punc\-tually) spectrum-regular if and only it is $0$-punc\-tually
cospectral, a concept which, as was claimed, can be defined through any of the
considered graph perturbations. In the next section, we generalize this result.

\section{$h$-Punctually walk-regular graphs}

To obtain some characterizations and properties of
$h$-punctually walk-regular graphs, we consider some basic
graph perturbations involving two vertices. With this aim, we
first perturb the vertices `separately', as done in the
previous section. Second, similar characterizations are derived
when we perturb the vertices `together'.

\subsection{Separate perturbations}

Let us first prove the following lemma concerning perturbations
{\bf P1-P3} for pairs of vertices in walk-regular graphs:
\begin{lema}
\label{co-sp} For all pairs of vertices $u,v$ and $w,z$ of a
walk-regular graph $G$, the following statements are
equivalent:
\begin{itemize}
\item[$(a)$]
 $\sp(G-u-v)=\sp(G-w-z)$.
\item[$(b)$]
 $\sp(G+uu+vv)=\sp(G+ww+zz)$.
\item[$(c)$]
 $\sp(G+u\ubar+v\vbar)=\sp(G+w\wbar+z\zbar)$.
\end{itemize}
\end{lema}
\proof The equivalence $(a)\iff (b)$ follows by applying
repeatedly  Eq. (\ref{G+uu vs G-u}) to obtain
$$
\phi_{G-u-v}(x)-\phi_{G+uu+vv}(x) = \phi_{G-u}(x)-\phi_{G+vv}(x),
$$
and using Proposition \ref{propo0}.
Analogously, from Eq. (\ref{G+uubar vs G-u}) we get
$$
\phi_{G-u-v}(x)-\phi_{G+u\ubar+v\vbar}(x) = x\phi_{G-u}(x)-x\phi_{G+v\vbar}(x),
$$
which proves $(a)\iff (c)$.
\endproof

Notice that, by this result and Proposition \ref{propo0},
each of the above conditions ($a$)-($c$) is equivalent to the
sets $\{u,v\}$ and $\{w,z\}$ being removal-cospectral. Moreover,
as mentioned before, this allows us to define $h$-punctually
cospectrality by requiring that every pair of vertices at
distance $h$ satisfies one of these conditions.

In turn, this leads to the following characterization of
$h$-punctually walk-regular graphs. It is, in a sense, a restatement of Corollary \ref{coroequiv}.

\begin{theo}
\label{theo-punct} For a walk-regular graph $G$ with diameter
$D$ and a given integer $h\le D$, the following statements
are equivalent:
\begin{itemize}
\item[$(a)$]
$G$ is $h$-punctually walk-regular.
\item[$(b)$]
$G$ is $h$-punctually spectrum-regular.
\item[$(c)$]
$G$ is $h$-punctually cospectral.
\end{itemize}
\end{theo}
\proof The equivalence $(a) \iff (b)$ was proved by the authors in
\cite{ddfgg10}. To prove the equivalence  $(b)\iff (c)$, we use Lemma
\ref{det2}, and follow the same line of reasoning as Rowlinson \cite{r96}.
Indeed, let $U=V\setminus\{u,v\}$ with $\dist(u,v)=h$, and $\C=x\I-\A$. Then,
\begin{eqnarray}
\phi_{G-u-v}(x) & = & \det \C[U] =\det \C\det \C^{-1}[\{u,v\}] \nonumber \\
 &  = & \phi_G(x) \det \left( \sum_{i=0}^d \frac{\E_i[\{u,v\}]}{x-\lambda_i}\right) \nonumber \\
 &  = & \phi_G(x) \det \left( \sum_{i=0}^d \frac{1}{x-\lambda_i}\left(
 \begin{array}{cc}
 m_{uu}(\lambda_i) & m_{uv}(\lambda_i) \\
 m_{uv}(\lambda_i) &  m_{vv}(\lambda_i)
 \end{array}\right)\right) \nonumber\\
 & = & \phi_G(x) \left[\left(\sum_{i=0}^d \frac{ m_{0i}}{x-\lambda_i}\right)^2
 -
 \left(\sum_{i=0}^d \frac{m_{uv}(\lambda_i)}{x-\lambda_i}\right)^2
 \right], \label{phi_G-{uv}}
\end{eqnarray}
where we have used that, as $G$ is walk-regular,
$m_{uu}(\lambda_i)=m_{vv}(\lambda_i)=m_{0i}$. Then, if $G$ is
$h$-punctually spectrum-regular, $m_{uv}(\lambda_i)=m_{hi}$
and, hence, $\phi_{G-u-v}(x)$ does not depend on $u,v$. This
proves $(b)\Rightarrow (c)$. Conversely, if
$\phi_{G-u-v}(x)=\phi_{G-w-z}(x)$ for some vertices $w,z$ at
distance $\dist(w,z)=h$, Eq. (\ref{phi_G-{uv}}) yields
$$
\left(\sum_{i=0}^d\frac{m_{uv}(\lambda_i)}{x-\lambda_i}\right)^2
=\left(\sum_{i=0}^d
\frac{m_{wz}(\lambda_i)}{x-\lambda_i}\right)^2
$$
for all $x \neq \lambda_0,\lambda_1,\ldots,\lambda_d$. Therefore,
$$
\sum_{i=0}^d\frac{m_{uv}(\lambda_i)}{x-\lambda_i}
=\pm \sum_{i=0}^d
\frac{m_{wz}(\lambda_i)}{x-\lambda_i}
$$
(since, as $p^2=q^2\Rightarrow p=\pm q$ holds for polynomials,
it also holds for rational functions). Consequently, taking
limits $x\rightarrow \lambda_i$, we have that either
$m_{wz}(\lambda_i)=m_{uv}(\lambda_i)$ for $i=0,1,\ldots,d$, or
$m_{wz}(\lambda_i)=-m_{uv}(\lambda_i)$ for $i=0,1,\ldots,d$.
But, since $m_{uv}(\lambda_0)=m_{wz}(\lambda_0)=\frac{1}{n}$,
we must rule out the second possibility and $G$ is
$h$-punctually spectrum-regular, thus proving that
$(c)\Rightarrow (b)$.
\endproof

\subsection{Joint perturbations}

We now consider the following perturbations involving two
given vertices $u,v$:
\begin{itemize}
\item[{\bf P4.}] $G\pm uv$ is the graph obtained from $G$
    by flipping the (non-)edge $uv$. (That is, changing the
    edge $uv$ into a non-edge or vice versa.)
\item[{\bf P5.}] $G_{u+v}$ is the (pseudo)graph obtained from $G$ by
    amalgamating the vertices $u$ and $v$ (if $u\sim v$ then the edge $uv$
    becomes a loop; if $u$ and $v$ have common neighbors, then multiple
    edges arise; the `new' vertex is denoted by $u+v$).
\item[{\bf P6.}] $G+u\ubar v$ is the graph obtained from
    $G$ by adding the 2-path $u\ubar v$ (thus creating a new
    so-called {\em bridging vertex} $\ubar$).
\end{itemize}
In the case that the graphs
$G+u\ubar v$ and $G+w\wbar z$ are cospectral, the pairs $(u,v)$
and $(w,z)$ are called {\em isospectral}; see Lowe
and Soto \cite{ls86}. In the following result, we show that for
walk-regular graphs, isospectral pairs can also be defined by
requiring cospectrality of the graphs obtained from perturbations {\bf
P4-P5}.

\begin{propo}
\label{iso-sp} Let $u,v$ and $w,z$ be pairs of vertices
of a walk-regular graph $G$ such that $u \sim v$ if and only if $w \sim z$.
Then the following statements are equivalent:
\begin{itemize}
\item[$(a)$]
  $\sp(G\pm uv)=\sp(G\pm wz)$.
\item[$(b)$]
  $\sp G_{u+v} =\sp G_{w+z}$.
\item[$(c)$]
  $\sp(G+u\ubar v)=\sp(G+w\wbar z)$.
\end{itemize}
\end{propo}
\proof We will prove that each of the above conditions is equivalent to
$m_{uv}(\lambda_i)=m_{wz}(\lambda_i)$, for all $i=0,1,\ldots,d$.
With respect to  $(a)$, note that, when $u\not\sim v$, the
adjacency matrix of the graph $G+uv$ can be written as
$$
\A_{+uv}=\left(
\begin{array}{ccc}
0 & 1 & \y^{\top} \\
1 & 0 & \z^{\top} \\
\y & \z & \A^*
\end{array}
\right)
$$
where $\A^*$ is the adjacency matrix of $G-u-v$. Then, by
applying twice Eq. (\ref{column-linearity}) (to the first
column and row) we have:
 \begin{eqnarray*}
 \det(\x\I-\A_{+uv}) & = & \det\left(
\begin{array}{ccc}
x & 0 & -\y^{\top} \\
0 & x & -\z^{\top} \\
-\y & -\z & x\I-\A^*
\end{array}
\right)+
\det\left(
\begin{array}{ccc}
0 & -1 & \vec0^{\top} \\
0 & x & -\z^{\top} \\
-\y & -\z & x\I-\A^*
\end{array}
\right) \\
 & + & \det\left(
\begin{array}{ccc}
0 & -1 & \vec0^{\top} \\
-1 & x & -\z^{\top} \\
\vec0 & -\z & x\I-\A^*
\end{array}
\right)+
\det\left(
\begin{array}{ccc}
0 & 0 & -\y^{\top} \\
-1 & x & -\z^{\top} \\
\vec0 & -\z & x\I-\A^*
\end{array}
\right).
\end{eqnarray*}
Thus,  with $\Psi_{uv}(x)$ denoting the $uv$-cofactor of
$x\I-\A$ (where $\A$ is the adjacency matrix of $G$), we get
\begin{equation}\label{phi(G+uv)}
\phi_{G+uv}(x) = \phi_G(x)-\phi_{G-u-v}(x)-2\Psi_{uv}(x).
\end{equation}
This equation was also derived by Rowlinson \cite{Row}. Moreover, using a
similar reasoning, Rowlinson \cite{r96} proved that, if $u\sim v$, we have
$$\phi_{G-uv}(x) = \phi_G(x)-\phi_{G-u-v}(x)+2\Psi_{uv}(x).$$
Then, from Eq. (\ref{phi_G-{uv}}) and since the $uv$-cofactor $\Psi_{uv}(x) $
can be computed as:
\begin{equation}\label{psiuv}
\Psi_{uv}(x)=\det(x\I-\A)((x\I-\A)^{-1})_{uv}=\phi_G(x)\sum_{i=0}^d \frac{m_{uv}(\lambda_i)}{x-\lambda_i},
\end{equation}
we get
\begin{equation}
\label{G+-uv}
\phi_{G\pm uv}(x)
=
 \phi_G(x)\left(1-\left(\sum_{i=0}^d \frac{ m_{0i}}{x-\lambda_i}\right)^2
 +
 \left(\sum_{i=0}^d \frac{m_{uv}(\lambda_i)}{x-\lambda_i}\right)^2
 \mp 2\sum_{i=0}^d \frac{m_{uv}(\lambda_i)}{x-\lambda_i}\right).
 \end{equation}
 Therefore $(a)$ is equivalent to
 $$
 \left(\sum_{i=0}^d \frac{m_{uv}(\lambda_i)}{x-\lambda_i} \mp 1 \right)^2=
 \left(\sum_{i=0}^d \frac{m_{wz}(\lambda_i)}{x-\lambda_i} \mp 1 \right)^2.
 $$
By the same reasoning as in the proof of Theorem
\ref{theo-punct}, this is equivalent to
$m_{uv}(\lambda_i)=m_{wz}(\lambda_i)$ for $i=0,1,\ldots,d$.

For $(b)$ we use similar techniques. Indeed, we now apply the formula
\begin{equation}
\label{G_{u+v}}
\phi_{G_{u+v}}(x)=\phi_{G-u}(x)+\phi_{G-v}(x)-(x-a_{uv})\phi_{G-u-v}(x)-2\Psi_{uv}(x),
\end{equation}
which is proved similarly as Eq. (\ref{phi(G+uv)}) (or by using $(2.8)$ and $(2.9)$ in Rowlinson \cite{r96} and Eq.
(\ref{G+uu vs G-u})). Using Eqs. (\ref{phi_G-{uv}}) and (\ref{psiuv}), it thus follows that
$(b)$ is equivalent to
$$
(x-a_{uv}) \left(\sum_{i=0}^d \frac{m_{uv}(\lambda_i)}{x-\lambda_i}\right)^2 -2\sum_{i=0}^d \frac{m_{uv}(\lambda_i)}{x-\lambda_i}=
(x-a_{uv}) \left(\sum_{i=0}^d \frac{m_{wz}(\lambda_i)}{x-\lambda_i}\right)^2 -2\sum_{i=0}^d \frac{m_{wz}(\lambda_i)}{x-\lambda_i},
$$
which, with $f(x)=\sum_{i=0}^d
\frac{m_{uv}(\lambda_i)}{x-\lambda_i}$ and $g(x)=\sum_{i=0}^d
\frac{m_{wz}(\lambda_i)}{x-\lambda_i}$, can be written as
$$
([x-a_{uv}][f(x)+g(x)]-2)[f(x)-g(x)]=0.
$$
The factor $[x-a_{uv}][f(x)+g(x)]-2$ in this equation cannot be zero. Indeed,
this could only happen if say $\lambda_j=a_{uv}$ and
$m_{uv}(\lambda_i)+m_{wz}(\lambda_i)$ would be 2 for $i=j$, and 0 otherwise.
This however leads to a contradiction by Lemma \ref{crossed-mul}(a). Hence,
$(b)$ is equivalent to $f(x)-g(x)=0$, which leads again to
$m_{uv}(\lambda_i)=m_{wz}(\lambda_i)$, $i=0,1,\ldots,d$.

Finally, case $(c)$ is managed by using the formula
$$
\phi_{G+u\ubar v}(x)= \phi_G(x) \left(x - 2\sum_{i=0}^d\frac{m_{0i}+m_{uv}(\lambda_i)}{x-\lambda_i} \right),
$$
(see Lowe and Soto \cite{ls86} or Rowlinson \cite[p. 216]{r96}).
\endproof

It is perhaps good to remind the reader that the condition
$m_{uv}(\lambda_i)=m_{wz}(\lambda_i)$ for all $i=0,1,\ldots,d$ implies that $u$ and $v$ are at the
same distance as $w$ and $z$ (by Lemma \ref{crossed-mul}$(c)$).
Inspired by this and the above result, we say that a graph $G$ with diameter
$D$ is {\em $h$-punctually isospectral}, for a given $h\le D$,
when every pair of vertices at distance $h$ satisfies one of
the conditions in Proposition \ref{iso-sp}. As a corollary of its proof, we then
obtain the following characterization of
$h$-punctually walk-regular (or $h$-punctually
spectrum-regular) graphs.

\begin{figure}[h!]
\begin{center}
\resizebox{40mm}{!}{\includegraphics{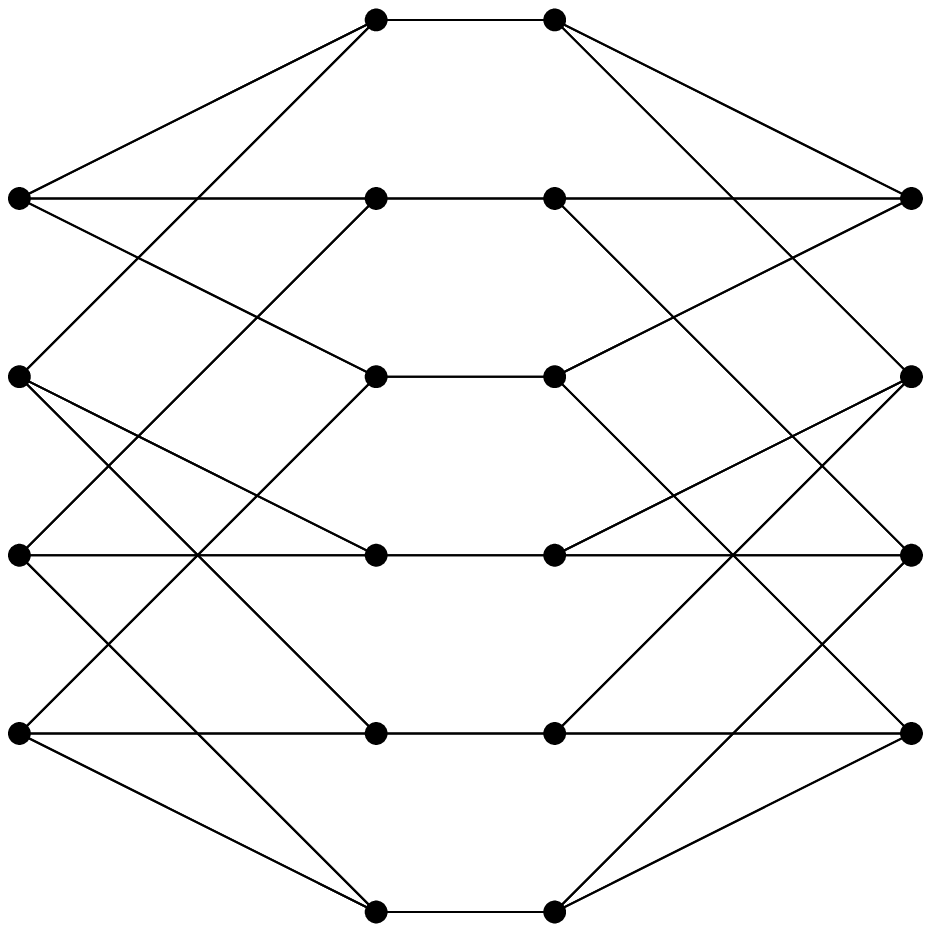}} \caption{Bussemaker-Cvetkovi\'c-Schwenk (twisted Desargues) graph} \label{twisted}
\end{center}
\end{figure}

\begin{figure}[h!]
\begin{center}
\resizebox{40mm}{!}{\includegraphics{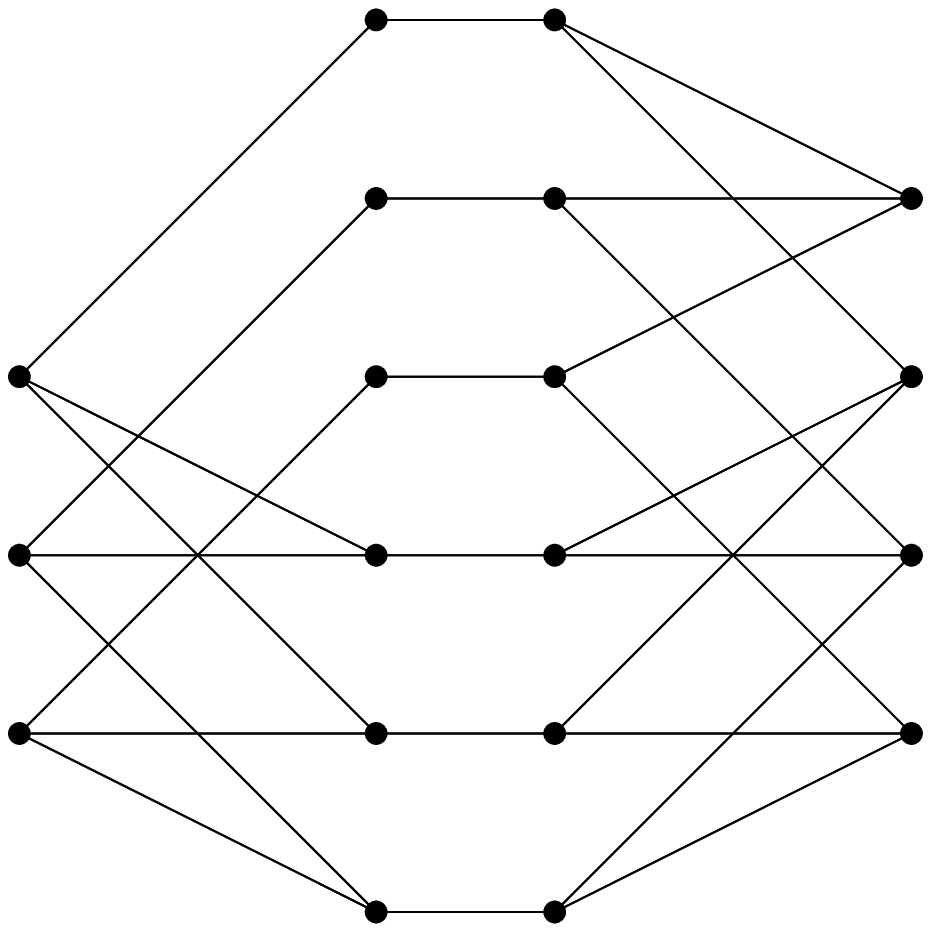}} \nobreak \hspace{-0.9 cm}
\resizebox{40mm}{!}{\includegraphics{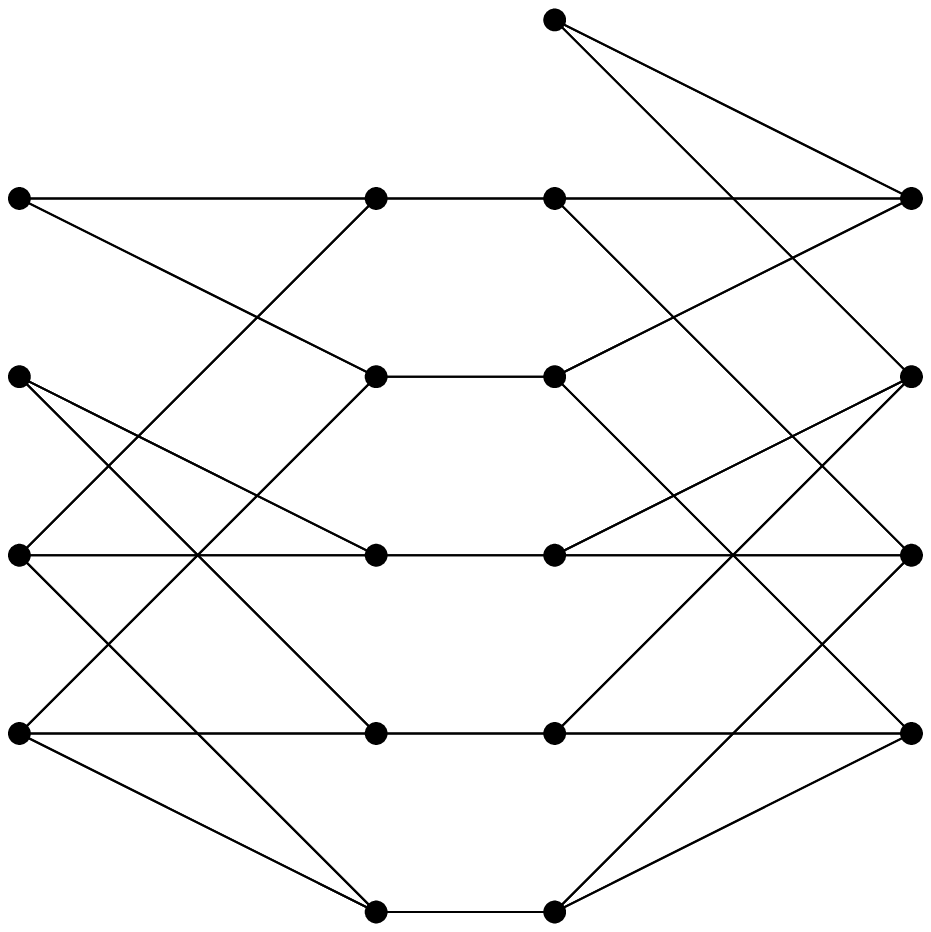}} \nobreak
\resizebox{40mm}{!}{\includegraphics{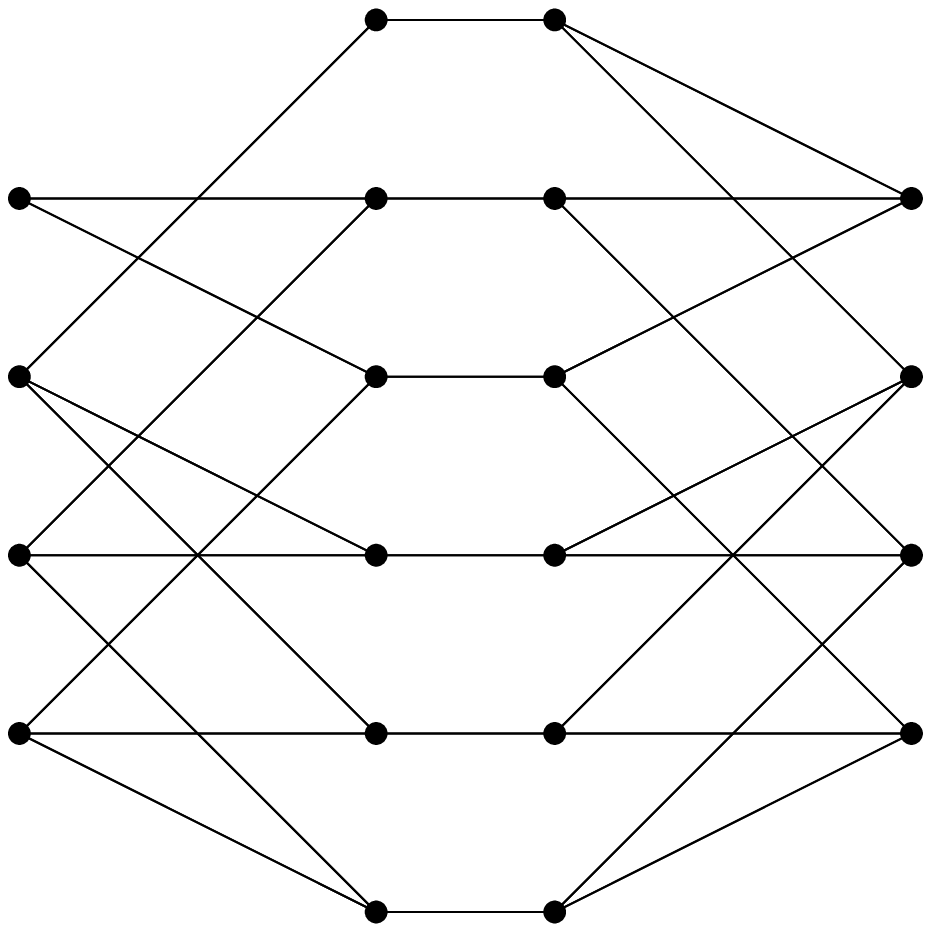}} \nobreak \hspace{-0.9 cm}
\resizebox{40mm}{!}{\includegraphics{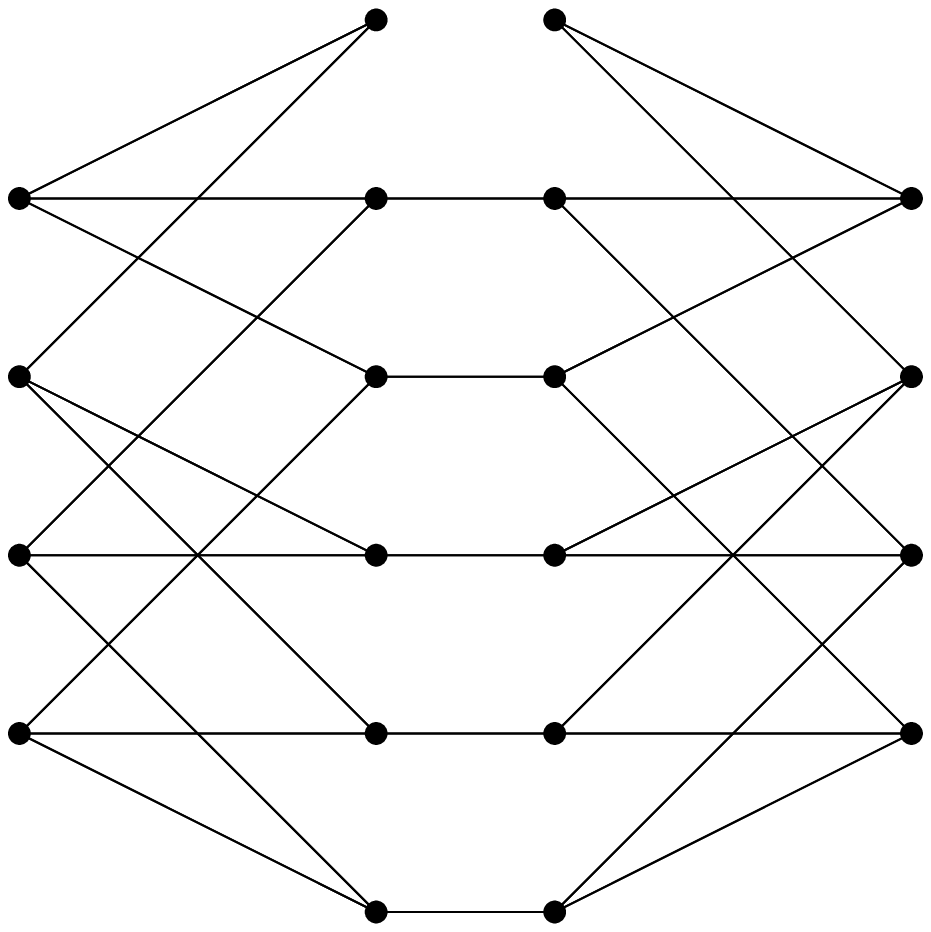}}
\end{center}
\caption{Pairs of cospectral graphs; removing vertices ($h=0$) and removing edges ($h=1$)}
\label{h0}
\end{figure}

\begin{figure}[h!]
\begin{center}
\resizebox{40mm}{!}{\includegraphics{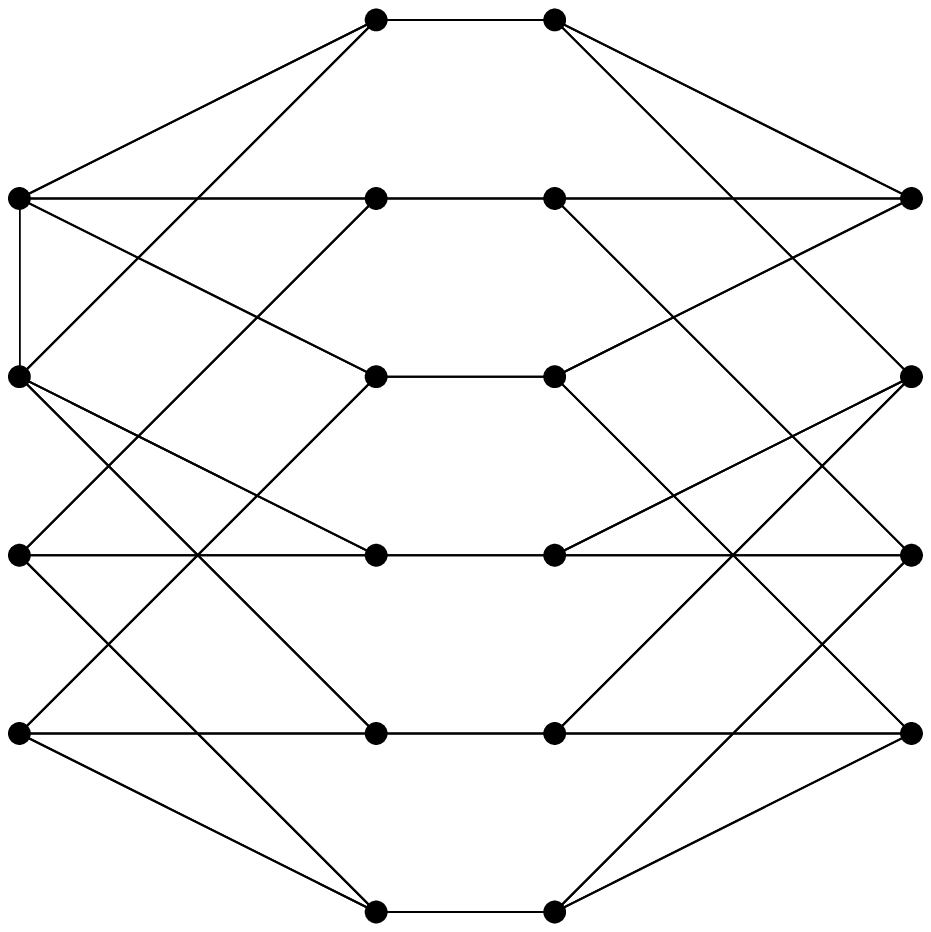}} \nobreak \hspace{-0.9 cm}
\resizebox{40mm}{!}{\includegraphics{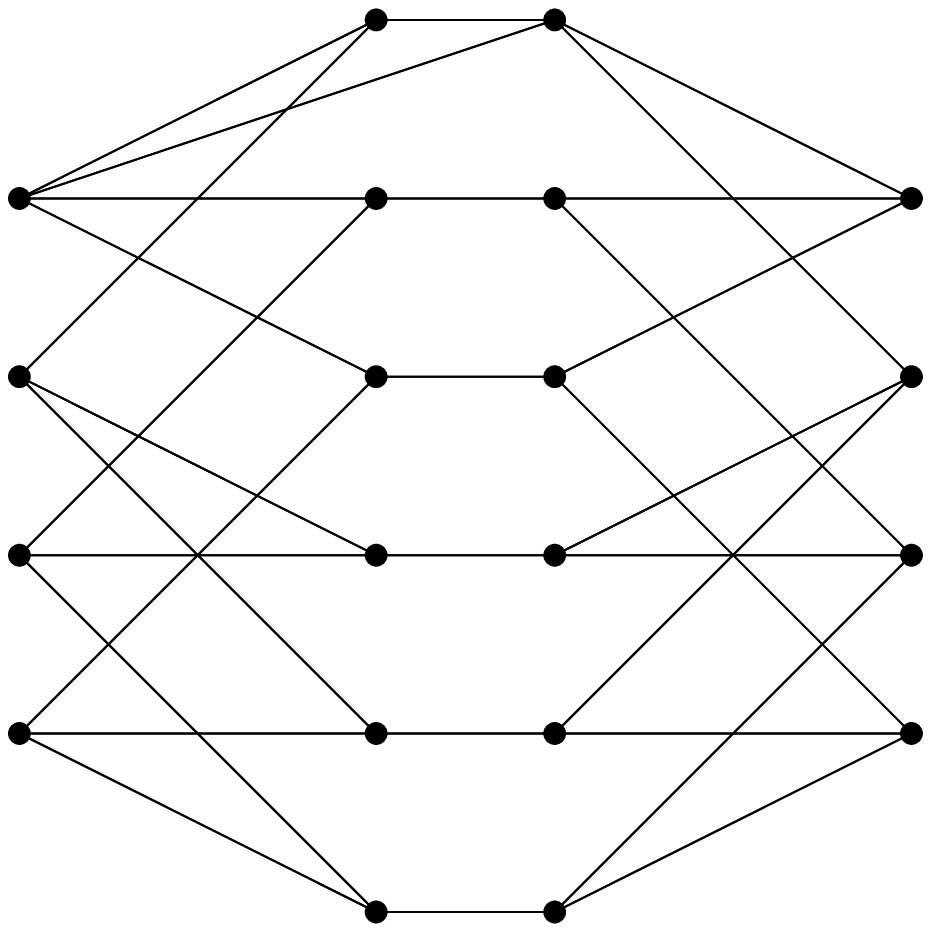}} \nobreak \hspace{-0.9 cm}
\resizebox{40mm}{!}{\includegraphics{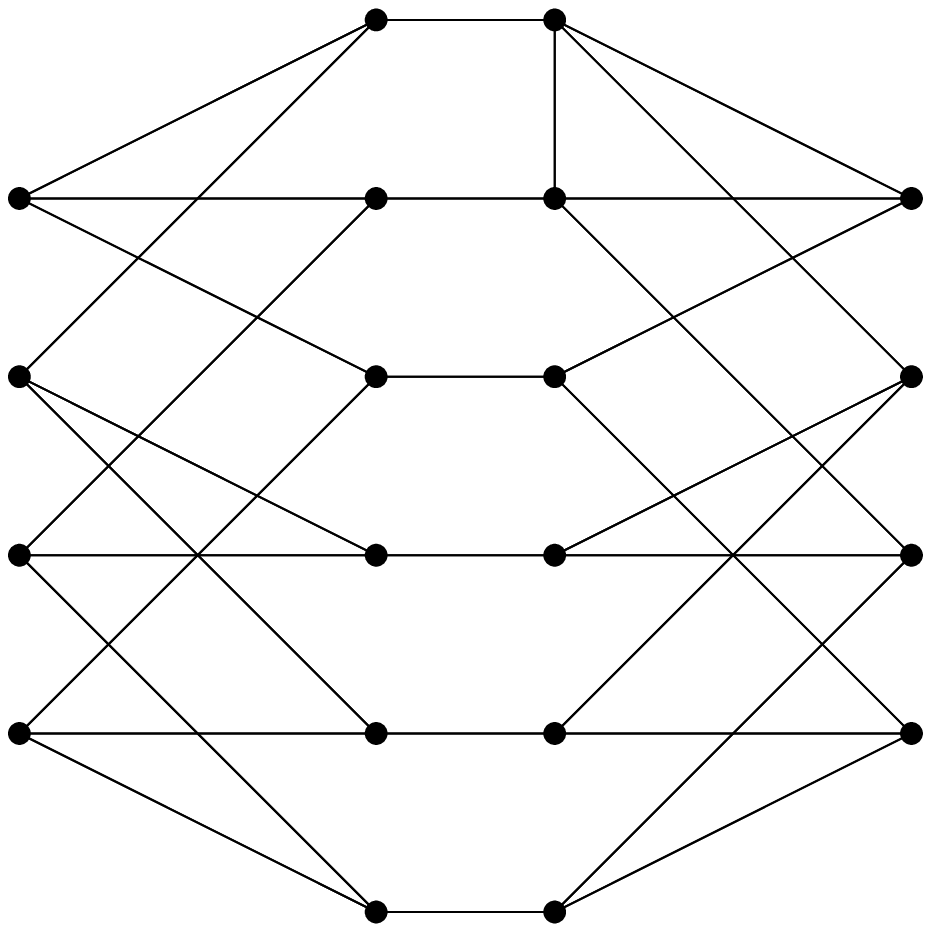}}
\end{center}
\caption{Triple of cospectral graphs; adding edges ($h=2$)}
\label{h2}
\end{figure}

\begin{figure}[h!]
\begin{center}
\resizebox{40mm}{!}{\includegraphics{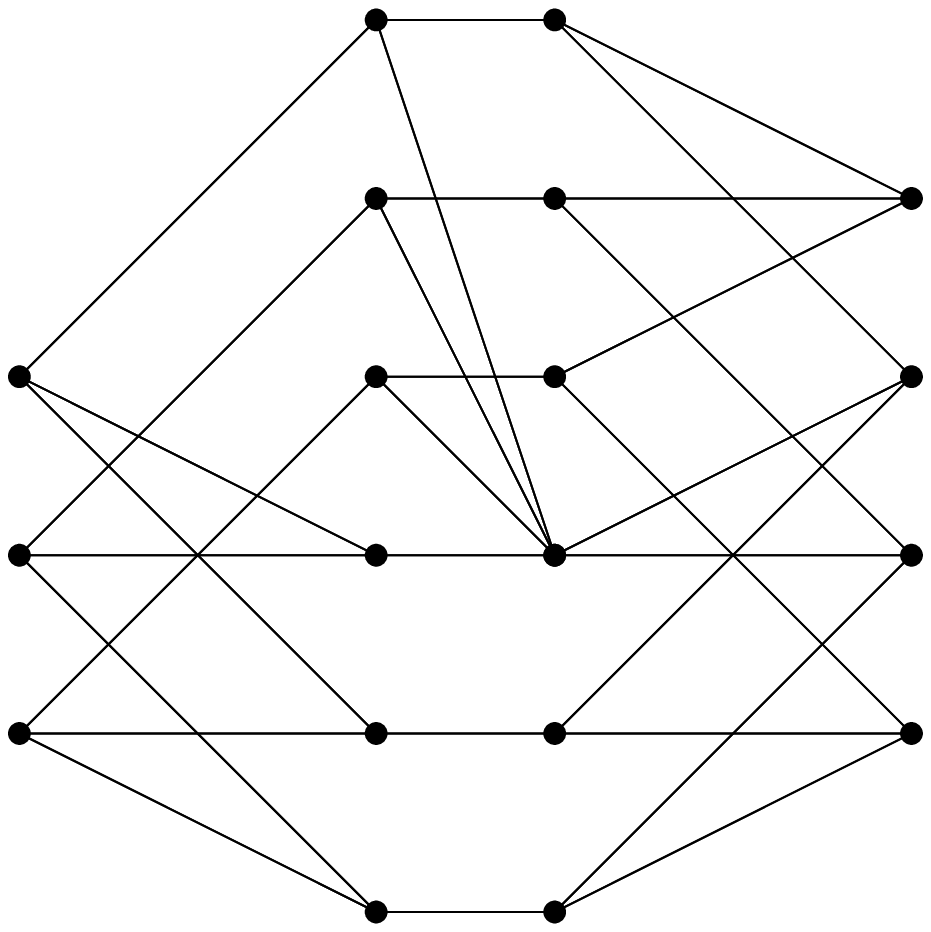}} \nobreak \hspace{-0.9 cm}
\resizebox{40mm}{!}{\includegraphics{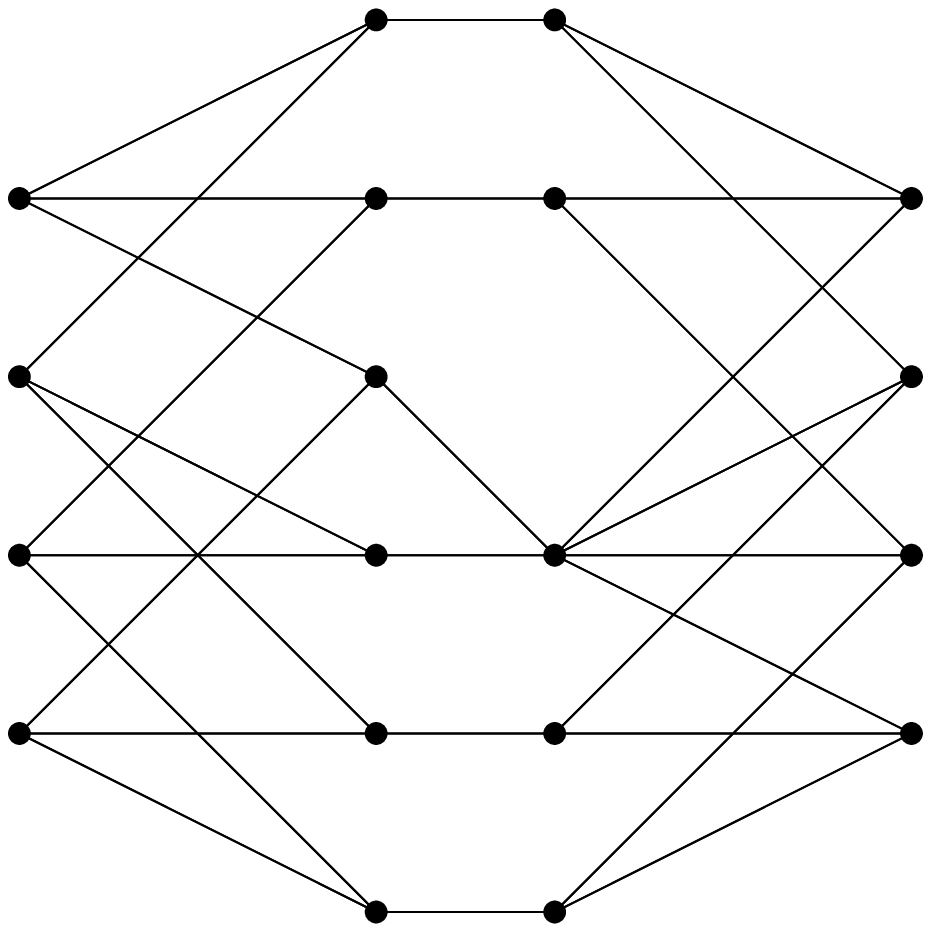}}
\end{center}
\caption{Pair of cospectral graphs; amalgamating vertices ($h=4$)}
\label{h4}
\end{figure}

\begin{coro}
\label{theo-punct-joint} For a walk-regular graph $G$ with
diameter $D$ and a given integer $h\le D$, the following
statements are equivalent:
\begin{itemize}
\item[$(a)$]
$G$ is $h$-punctually walk-regular.
\item[$(b)$]
$G$ is $h$-punctually spectrum-regular.
\item[$(c)$]
$G$ is $h$-punctually isospectral.
\end{itemize}
\end{coro}
We finish this section with an example of an almost
distance-regular graph that can be used to produce many kinds
of cospectral graphs by applying the above perturbations. The
graph we use is one of the thirteen cubic graphs with integral
spectrum. These graphs were classified by Bussemaker and
Cvetkovi\'c \cite{BC} and Schwenk \cite{Sch}. Of these graphs
we take the one that is cospectral (with spectrum
$\{\pm3^1,\pm2^4,\pm1^5\}$), but not isomorphic, with the
Desargues graph; see Figure \ref{twisted}. This graph can be
obtained by switching from the Desargues graph (take the four
right-most vertices as switching set) and also by twisting it
(in a similar way as in the twisted Grassmann graphs of
\cite{DK}); cf. \cite[Sect. 3.2]{vdhks06}. It is a bipartite
graph with diameter $D=5$ that is almost distance-regular in
the sense that it is $h$-punctually walk-regular for all $h$
except $h=3$. It has two orbits of vertices under the action of
its automorphism group; the middle twelve vertices are
different from the others. This means that if we remove a
vertex from the middle, and remove a vertex from the left four,
we obtain non-isomorphic graphs that are cospectral (apply the
above with $h=0$). These graphs are shown as the left pair in
Figure \ref{h0}. There are two kinds of edges, three kinds of
pairs of vertices at distance 2, and two kinds of vertices at
distance 4. These (for example) give cospectral graphs as shown
on the right in Figure \ref{h0}, and in Figures \ref{h2} and
\ref{h4}, respectively. There is only one kind of pair of
vertices at distance 5, so these cannot be used to get
non-isomorphic but cospectral graphs. We finally remark that
this example can be generalized easily to other twisted graphs
that are described in \cite[Sect. 3.1-2]{vdhks06}; for example
the distance-regular twisted Grassmann graphs.

\subsection{Multiple perturbations}

For the sake of simplicity, we have only considered perturbations in a single graph $G$ so far.
One could however also use the above perturbations in cospectral graphs $G$ and $G'$ to get new cospectral graphs (as is well known from the literature).
The conditions for this to work are similar as before: the crossed local multiplicities
$m_{uv}(\lambda_i)$ (in $G$) and $m'_{wz}(\lambda_i)$ (in $G'$) should be the same for all $i=0,1,\ldots,d$ (or alternatively: the number of walks $a_{uv}^{(\ell)}$ (in $G$) and $a_{wz}'^{(\ell)}$ (in $G'$) should be the same for all $\ell$).

Consider the Desargues graph, for example. This graph is cospectral to the above mentioned twisted Desargues graph.
By removing a vertex ($h=0$), removing an edge ($h=1$), adding an edge $(h=2)$, and amalgamating vertices ($h=4)$ in the Desargues graph, one gets cospectral graphs of the graphs in the above figures. One could even exploit the case $h=5$ now.

For the next step --- multiple perturbations --- it is hard to avoid working
with different (but cospectral) graphs. We next consider removal-cospectral
sets $U,U'$ belonging to cospectral (but not necessarily isomorphic) graphs
$G,G'$ (i.e., there exists a one-to-one mapping $U \rightarrow U'$ such that,
for every $W\subset U$, the graphs $G-W$ and $G'-W'$ are cospectral), as is
usually done in the literature. The following proposition shows that all
perturbations {\bf P1-P6} leave the property of two sets being
removal-cospectral invariant, and gives new insight into some of the previous
implications.

\begin{propo}
Let $U$ and $U'$ be removal-cospectral sets in cospectral graphs $G$ and $G'$,
and let $u,v \in U$ with corresponding vertices $u',v' \in U'$. Let
$\widetilde{U},\widetilde{U}'$ be the sets obtained from $U,U'$ after
perturbing vertices $u$ and $u'$ according to one of the perturbations {\bf
P1-P3}, or perturbing pairs of vertices $u,v$ and $u',v'$ through one of the
perturbations {\bf P4-P6}, where possible new vertices $u+v$, $\ubar$, $\ubar'$
are included in $\widetilde{U},\widetilde{U}'$. Let $\widetilde{G}$ and
$\widetilde{G}'$  be the resulting perturbed graphs. Then, the sets
$\widetilde{U},\widetilde{U}'$ are removal-cospectral in $\widetilde{G}$ and
$\widetilde{G}'$.
\end{propo}
\proof We only prove the result for amalgamation (that is, {\bf
P5}), as the other cases are either very simple, or similar, or follow from
Schwenk's results in \cite{sch79}. Thus, let us amalgamate
$u,v\in U$ and $u',v'\in U'$ to obtain  $\widetilde{G}=G_{u+v}$ and
$\widetilde{G}'=G'_{u'+v'}$. Now, consider a subset $S\subset
\widetilde{U}$ and its corresponding set $S'\subset \widetilde{U}'$. We should prove that
$G_{u+v}-S$ and $G'_{u'+v'}-S$ are cospectral. To do this, we
must consider two cases: If $u+v\in S$, then $G_{u+v}-S =
G-(S\cup\{u,v\})$ and $G'_{u'+v'}-S' = G'-(S'\cup\{u',v'\})$.
Hence, these two graphs are cospectral. Otherwise, if
$u+v\not\in S$, then $G_{u+v}-S=(G-S)_{u+v}$ and $G'_{u'+v'}-S'
= (G-S')_{u'+v'}$, and these graphs are also cospectral because
$U\setminus S$ and $U'\setminus S'$ are removal-cospectral in
$G-S$ and $G'-S'$ (notice that, since $u,v\in U\setminus S$ and
$u',v'\in U'\setminus S'$, we can apply Eq. (\ref{G_{u+v}}) or
repeat the above argument).
\endproof

As a consequence, notice that the different one-vertex and two-vertex
perturbations can be repeated over and over again to obtain different
cospectral graphs $\widetilde{G}$ and $\widetilde{G}'$. In other words, from
two removal-cospectral sets $U,U'$, one can, for example, amalgamate several
vertices, or combine amalgamation with other operations such an edge
removal/addition (hence also contract an edge), adding pendant edges, etc., to
obtain new removal-cospectral sets $\widetilde{U},\widetilde{U}'$ in the
corresponding cospectral graphs $\widetilde{G},\widetilde{G}'$. This suggests
the following definition: Two vertex subsets $U,U'$ of cospectral graphs $G,G'$
are called {\em perturb-cospectral} if for all subsets $S\subset U$ and
$S'\subset U'$, the perturbed graphs $\widetilde{G}$ and $\widetilde{G}'$,
obtained by applying {\bf P1-P6} to corresponding vertices of $U$ and $U'$, are
cospectral.

\section{Distance-regular graphs}

In this section we use the above results to obtain some new
characterizations of distance-regular graphs.

In \cite{ddfgg10}, the authors considered also the following
concepts: A graph $G$ is {\em $m$-walk-regular} (respectively
{\em $m$-spectrum-regular}) when it is $i$-punctually
walk-regular (respectively $i$-punctually spectrum-regular) for
every $i\le m$. Similarly, we say that $G$ is {\em
$m$-cospectral} (respectively, {\em $m$-isospectral}) when it
is $i$-punctually cospectral (respectively, {\em $i$-punctually
isospectral}) for every $i\le m$. Using these definitions,
Theorem \ref{theo-punct} and Corollary \ref{theo-punct-joint} have the
following direct consequence:
\begin{coro}
\label{coro-part} For a walk-regular graph $G$ with diameter
$D$ and a given integer $m\le D$, the following statements are
equivalent:
\begin{itemize}
\item[$(a)$]
$G$ is $m$-walk-regular.
\item[$(b)$]
$G$ is $m$-spectrum-regular.
\item[$(c)$]
$G$ is $m$-cospectral.
\item[$(d)$]
$G$ is $m$-isospectral.
\end{itemize}
\end{coro}
Moreover, as mentioned in Section \ref{secgraphspectra},
Rowlinson \cite{r97} proved that a graph $G$ is
distance-regular if and only if it is $D$-walk-regular. Hence,
we get the following characterization:
\begin{theo}\label{drgDcospectral}
Let $G$ be a graph with diameter $D$. Then, the following
statements are equivalent:
\begin{itemize}
\item[$(a)$]
$G$ is distance-regular.
\item[$(b)$]
$G$ is $D$-cospectral.
\item[$(c)$]
$G$ is $D$-isospectral.
\end{itemize}
\end{theo}
In fact, notice that we also proved the following result:
\begin{theo}
A graph $G=(V,E)$ is distance-regular if and only if every two
isometric subsets  $U,U'\subset V$ are perturb-cospectral.
\end{theo}
Part of the case $D=2$ of Theorem \ref{drgDcospectral} was already observed by
Cvetkovi\'c and Rowlinson \cite{CR}; they showed that if $G$ is strongly
regular, then $\phi_{G-u-v}(x)$ depends only on whether or not $u$ and $v$ are
adjacent. See also the observation by Godsil on cospectral graphs in strongly
regular graphs in \cite[Prop. 8]{DH}.\\

\noindent {\bf Acknowledgement} The authors are grateful to Ernest Garriga,
Willem Haemers, and Peter Rowlinson for discussions on the topic of this paper.
They also thank an anonymous referee for several useful comments. Research
supported by the Ministerio de Educaci\'on y Ciencia, Spain, and the European
Regional Development Fund under project MTM2008-06620-C03-01 and by the Catalan
Research Council under project 2009SGR1387.

\end{document}